# The Fractal Nature of Maps and Mapping

Bin Jiang

Department of Technology and Built Environment, Division of Geomatics
University of Gävle, SE-801 76 Gävle, Sweden
Email: bin.jiang@hig.se

**Abstract**
A fractal can be simply understood as a set or pattern in which there are far more small things than large ones, e.g., far more small geographic features than large ones on the earth surface, or far more large-scale maps than small-scale maps for a geographic region. This paper attempts to argue and provide evidence for the fractal nature of maps and mapping. It is the underlying fractal structure of geographic features, either natural or human-made, that make reality mappable, large-scale maps generalizable, and cities imageable. The fractal nature is also what underlies the beauty of maps. After introducing some key fractal concepts such as recursion, self-similarity, scaling ratio, and scaling exponent, this paper demonstrates that fractal thought is rooted in long-standing map-making practices such as series maps subdivision, visual hierarchy, and Töpfer's radical law. Drawing on previous studies on head/tail breaks, mapping can be considered a head/tail breaks process; that is to divide things around an average, according to their geometric, topological and/or semantic properties, into the head (for those above the average) and the tail (for those below the average), and recursively continue the dividing process for the head for map generalization, statistical mapping, and cognitive mapping. Given the fractal nature of maps and mapping, cartography should be considered a perfect combination of science and art, and scaling must be formulated as a law of cartography or that of geography in general.

**1. Introduction**
The word 'fractal' refers to things that are fragmented, irregular, and not smooth (Mandelbrot 1982). In fact, the term literally sounds very much like 'fraction'. Both 'fractal' and 'fraction' come from the same Latin adjective *fractus*, meaning 'broken', but there are some fundamental differences between the two. For the former, a whole is broken into a large number of different irregular pieces, whereas for the latter, a whole is broken into a small number of equal, regular pieces. Given the definition of 'fractal' and 'fraction', the things surrounding us can be divided into two basic types: those with regular and smooth shapes that can be measured by Euclidean geometry, and those with irregular shapes that Euclidean geometry fails to characterize. Let us do a simple experiment. A piece of glass of regular shape (such as a circle or rectangle) is clearly measurable in its entirety. Now, throw the glass onto stone ground, and it is very likely to be broken into a considerable number of irregular pieces. Unlike the piece of glass itself, the broken pieces are difficult to measure accurately, due to their irregular shapes. For the sake of simplicity, we measure their maximum lengths (in any direction) to provide information on their sizes. Eventually, after measuring all of the broken pieces, we realize that there are far more small pieces than large ones. This paper intends to argue that maps share the same scaling or fractal property as the broken pieces, and mapping is essentially to reflect or reveal the underlying fractal structure.

A map is a visual representation of what is too large and/or too complex to be perceived in its entirety. In the present context, the too large or too complex, or reality in general, widely refers to the earth, the moon, and any celestial bodies (which are too large), as well as virtual spaces such as the Internet, the World Wide Web, and social media (whose structures are too complex) (e.g., Jiang and Ormeling 1997, Dodge and Kitchin 2001). Mapping involves map-making processes such as generalization, classification, and symbolization so as to reduce and simplify reality into a two- or three-dimensional space. It should be noted that maps are also used for what is too small, like the human brain (which is also too complex), in which case magnified map scales should be adopted (Carter 2010). Such things as mountains, rivers, roads, and settlements are essentially visible in reality, and they constitute the major content of topographic maps (Robinson et al. 1995). On the other hand, many natural and human-made things are invisible, such as precipitation, temperature, the Internet, and population density, and these become visible through thematic or statistical mapping (Slocum et al. 2008, Kraak and Ormeling 2010). Things to be represented in a map are not decided arbitrarily; they are determined by several factors such as map scale, map theme, and map usage, with respect to their geometric (location, size and direction), topological (connectivity or popularity) and semantic (historical or societal meaning) properties.

Over the past 100 years, maps and mapping have been scientifically studied from various perspectives such as map design, map use, and even map perception. This can be seen in works by some of the preeminent cartographers and thinkers, which include: *'On the Nature of Maps and Map Logic'* (Eckert 1908), *'The Nature of Maps'* (Robinson and Petchenik 1976), *'Semiology of Graphics'* (Bertin 1983), and *'How Maps Work'* (MacEachren 2004). This is a non-exclusive list, for there are many related works in the literature that focus on scientific understanding of maps and mapping. In the same spirit, this paper attempts to argue and provide evidence for the fractal nature of maps and mapping. It is the underlying fractal structure of geographic features, or equivalently the presence of the far more small things than large ones in general, that makes reality mappable, large-scale maps generalizable, and cities imageable. The fractal nature is also what underlies the beauty of maps. It is important to note that the notion of far more small things than large differs fundamentally from that of more small things than large ones: the former represents a nonlinear relationship, while the latter a linear relationship. Note also that the term 'scale' is used in the context of this paper with completely opposite meanings: the map scales (the ratio of a distance on the map to the corresponding distance on the ground) on the one hand and the measuring scales on the other. A large-scale map implies a small measuring scale, but a small-scale map means a large measuring scale; assuming that the minimum visual resolution of a map is 1 millimeter, the measuring scales (or yardsticks) for three map scales (1:250K, 1:500K, and 1:1M) are respectively 250 meters, 500 meters, and 1,000 meters. In this paper, 'scale' is primarily used with reference to 'measuring scale', except where it appears alongside 'map'.

The remainder of this paper is structured as follows. Section 2 introduces some key fractal concepts using classic fractals such as the Koch curve and Fibonacci numbers. Section 3 examines how fractal thought is rooted in map-making practices such as series maps subdivision, visual hierarchy, and Töpfer's radical law. Section 4 illustrates how mapping in general, or map generalization, statistical mapping, and cognitive mapping in particular, can be considered a head/tail breaks process. Section 5 further discusses the implications of the fractal nature of maps and mapping for cartography and geography in the era of big data. Finally Section 6 concludes the paper.

**2. Fractal concepts: recursion, self-similarity, scaling ratio, and scaling exponent**
All maps use reduced scales. In fact, there is no need for 1:1 scale maps, since the power of maps lies in reduced map scales. Maps capture the same structure of the territories they represent, and allow us to see the territories in their entirety while ignoring trivial things. Importantly, different scales of maps constitute a recursive relationship, for example, a map of 1:1M is the map of four maps of 1:500K, and 16 maps of 1:250K. This recursive, or nested, relationship is key to understanding the fractal nature of maps, and the recursion is also closely related to other fractal concepts such as self-similarity, scaling ratio (or equivalently, similarity ratio), and scaling exponent (or fractal dimension). Here, we use some classic fractals such as the Koch curve and Fibonacci sequence, or the golden rectangles, as working examples to illustrate these rather abstract concepts.

In 1904, the Swedish mathematician Helge von Koch (1870-1924) invented what is now called the Koch curve (Figure 1). Given a straight line of one unit (called an initiator), *divide* it equally into three parts, and *replace* the middle part with two sides of an equilateral triangle, resulting in a curve made of four line segments, each of which is one third of the initiator in length (the resulting four segments are collectively called a generator). This is the first iteration. The result of the first iteration, the curve of four line segments, is fed back to form the input for the second iteration, and each of the four segments then repeats the same dividing and replacement processes as in the first iteration. The second iteration leads to a curve made of 16 line segments, each of which is one ninth of the initiator in length. Note that this is a recursive process, which implies that the result of the previous iteration is fed back to form the input of the subsequent iteration (see the left panel of Figure 1 for the first three iterations). Importantly, the first-iteration curve is embedded in the second-iteration curve, and the first two curves are embedded in the third-iteration curve, thus forming the cascade structure of the Koch curve, of the Koch curve, of the Koch curve, and so on endlessly. For example, the curve of the third iteration includes 4, 16, and 64 line segments with respect to the three scales: 1/3, 1/9, and 1/27. The so-called Koch curve (singular) refers to the curve when the line segment approaches an infinitely small size, so that the corresponding curve has an infinite length. The Koch curve is a good example of recursion. Another good example is the Fibonacci sequence: 1, 1, 2, 3, 5, 8, 13, 21, 34, 55, 89, ..., in which each subsequent number (except for the first two) is the sum of the previous two. The Fibonacci sequence is the foundation of the golden rectangles, or the golden ratio phi; refer to the right panel of Figure 1 for the first four iterations.

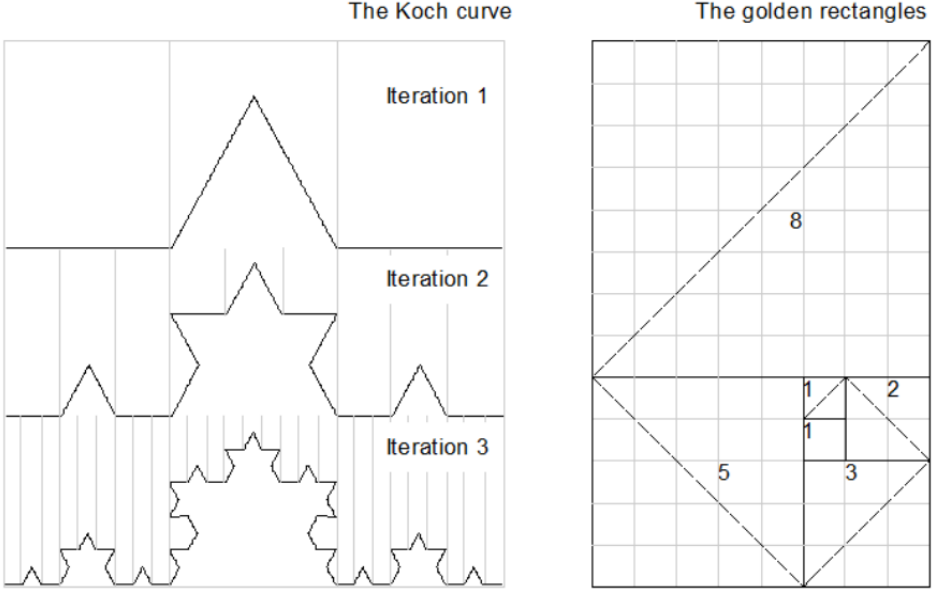

Figure 1: The Koch curve (left) and the golden rectangles (right) from the Fibonacci sequence (Note: Start with a line of one unit, divide it equally into three thirds, and replace the middle one by two sides of an equilateral triangle. The dividing and replacing process goes on recursively, ending up with the Koch curve, named after the Swedish mathematician Helge von Koch (1870-1924). The golden rectangles represent part of the Fibonacci sequence: 1, 1, 2, 3, 5, and 8.)

The Fibonacci sequence has a very interesting property, called self-similarity. The ratio of the subsequent number to the previous number is respectively 1, 2, 1.5, 1.667, 1.600, 1.625, 1.615, 1.619, 1.618, 1.618, …., which is persistently approaching to the golden ratio phi = 1.618. This self-similarity property is clearly reflected in the golden rectangles, in which all rectangles are self-similar to the largest one or the whole, as they tend to have the same length-to-width ratio. The ratio also has a more general name, which is called 'scaling ratio'. The reader might have noticed that in the course of generating the Koch curve, the line segment is decreased by one third. This one third is the scaling ratio, and it remains the same for all of the iterations. Every part of the Koch curve is self-similar to the whole curve. Note that the part is not defined arbitrarily. If one part is not self-similar to the whole, then the part is not rightly defined. This sounds circular, but it makes perfect sense with respect to the Koch curve. Put in more general terms, self-similarity indicates that a part has the same shape as the whole.

As mentioned earlier, the third iteration of the Koch curve consists of 4, 16, and 64 line segments with respect to the three scales 1/3, 1/9, and 1/27. Note that we must take a recursive perspective; otherwise, one would only see 64 line segments. When the two datasets (the number of line segments and the scales) are transferred into an Excel sheet, and plotted as a scatter chart, a power law relationship emerges between the two datasets, with the power law exponent 1.262. This chart is called the Richardson plot (Richardson 1961). The power law exponent is also called the 'scaling exponent', or 'fractal dimension'. The scaling exponent is a highly important concept in fractal geometry, and can be compared to the concept of scale in Euclidean geometry. Parameters such as radius and side length are sufficient to characterize regular shapes, such as circles and rectangles. For irregular shapes like the Koch curve, or fractals in general, the length or size is not measurable, and depends on the scale of the measuring yardstick, which is known as the 'conundrum of length' or the 'Steinhaus Paradox' (Richardson 1961, Perkal 1966) (c.f., Section 3 for more details). Because of this, the Koch curve had been called 'monstrous' or 'pathological', before fractal geometry was developed (Mandelbrot 1982).

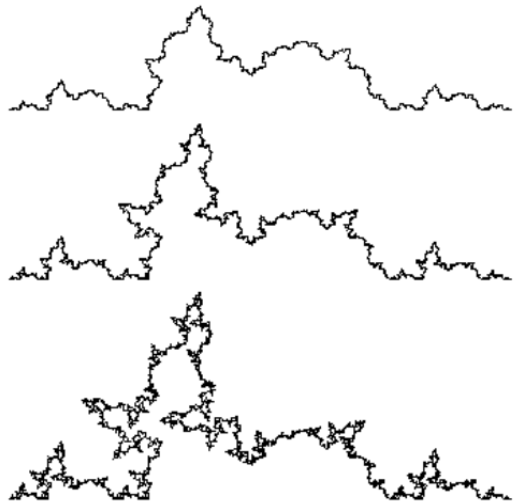

Figure 2: Three randomized or statistical Koch curves
(Note: The curves look very much like clouds or convoluted coastlines. The scaling exponent of the curves decreases from bottom to top)

The aforementioned fractal concepts, recursion, scaling ratio, self-similarity, and scaling exponent, are closely inter-related. The recursion helps to explain the generation process of the Koch curve; the scaling ratio helps to generate smaller structures that are self-similar to the whole; and the scaling exponent is of use for characterizing the complexity of the Koch curve, or fractals in general. For the sake of simplicity, we introduce these concepts using the classic fractals. In fact, none of the classic fractals appear in nature or society, in particular with respect to the infinite scaling range (either infinitely small or infinitely large). Fractals seen in nature and society (1) are not strictly self-similar, but statistically self-similar; and (2) do not fall within an infinite scaling range, but within a limited range. This insight into fractals was first developed by Benoit Mandelbrot (1967) in his classic work *'How Long Is the Coast of Britain?'*. Fractals in nature and society are more like coastlines rather than Koch curves. In other words, the appearance of Koch curves is not natural enough. By introducing some randomness, Koch curves look very much like coastlines or clouds (Figure 2). The change from strict self-similarity to statistical self-similarity, or equivalently from classic fractals to statistical fractals, is not trivial, and it took about 100 years for the new field of fractals to come into being. Now with fractal geometry, we can see that the Koch curve and other classic fractals are neither 'monstrous' nor 'pathological', since there is a hidden order beneath the surface complexity of fractals. Their complexity can be characterized by the scaling exponent.

## 3. Fractal thought in cartography

Throughout the thousands of years of map-making history, cartographers have been *unconsciously* guided by fractal thought in map design, map reading, and map related research such as map generalization. Here, we have deliberately used the adjective *'unconsciously'*, since fractal thought was developed prior to fractal geometry, and was not directly influenced or inspired by fractal geometry. In what follows, we will review some cartographic principles that manifest fractal thought such as topographic series maps subdivision, visual hierarchy in map design and reading, and Töpfer's radical law for map generalization. We begin with what is known as 'Steinhaus Paradox', the fact that the measured length of geographic features such as rivers increases with increasing map scales, or equivalently decreases with measuring scales (or yardsticks) (Richardson 1961, Steinhaus 1960, cited in Goodchild and Mark 1987). In the 1960s, the 'Steinhaus Paradox' received much attention in the fields of cartography and geography in relation to understanding the paradox and for measuring things in maps (e.g., Perkal 1966, Nystuen 1966). The paradox or myth was finally uncovered following publication of *'How Long is the Coast of Britain?'*, and subsequently, the establishment of the new field of fractal geometry (Mandelbrot 1967, 1982).

Topographic series maps (of different scales) are subdivided according to certain constant scaling ratios. If a geographic region is covered by one 1:1M map, then there would be four 1:500K maps, and 16 1:250K maps for the region. The scaling ratio of the series maps 1:250K, 1:500K, and 1:1M is 1/2. Similarly, a 1:450K map is derived from nine 1:150K maps, or from 81 1:50K maps, and the scaling ratio is 1/3. The series maps 1:50K, 1:150K, and 1:450K can be compared to the different scales of 1/3, 1/9 and 1/27 Koch curves, as both have the same scaling ratio of 1/3. There are far more large-scale maps than small-scale maps; this sounds obvious, but it in fact manifests some naïve fractal thought. The number of map sheets increases exponentially as the map scale decreases, just as that the number of boxes increases exponentially as the box size decreases with the box-counting method for calculating fractal dimension. From the fractal perspective, map scales among a series of maps must maintain a constant scaling ratio, which involves an integer like 2 or 3. In theory, map scales for series maps can be arbitrarily defined between the largest and the smallest, but they create a potential problem for map generalization; further discussion on this is provided below and in Section 5.

Things on a map are not all displayed at the same visual layer, and they tend to be organized hierarchically through color, line thickness, and size of symbols. This is referred to as 'visual hierarchy', which is a very important cartographic principle, not only for map design, but for map reading as well (Robinson et al. 1995, Kraak and Ormeling 2010, Slocum et al. 2008). For example, Google Maps uses yellow, white, and grey mixed with the line thickness to show the three hierarchical levels of street at the city scale. Visual hierarchy reflects the fact that there are far more small things than large ones on a map, and thus it embodies some fractal thought. Most important things (which are usually very few) should appear at the highest visual level, followed by those that are less important, with the vast majority of trivial things at the ground level. Visual hierarchy enables us to read a map as a whole rather than as individually disconnected elements, forms, and shapes; the whole is more than the sum of its parts, as advocated by Gestalt psychology (Köhler 1992). Map design and map reading should try to avoid the kind of ambiguity in figure-ground perception, famously illustrated by the 'Rubin vase' (Rubin 1921). One effective way in which to avoid such ambiguity is to adopt head/tail breaks (Jiang 2013a) in order to illustrate the underlying hierarchy; refer to the following section for more details on head/tail breaks.

Töpfer's radical law is an empirical law, also widely known as the Principle of Selection, that governs the number of map objects to be selected from a large-scale map (called the source map) to a small-scale map (called the derived map), e.g,, from 1:500K to 1:1M (Töpfer and Pillewizer 1966). This law states that the ratio of the number of objects to be selected to the number of objects in the source map is not simply inversely proportional to the two corresponding map scales. For example, from 1:500K to 1:1M, the selected objects are not just 1/2 of the objects in the source map, but rather calculated according to the square root of 1/2. Different versions of Töpfer's radical law can be applied, depending on other factors. This nonlinear relationship between the selected number and the source number constitutes the essence of Töpfer's radical law. This nonlinearity reflects some fractal thought. However, the Principle of Selection was initially obtained empirically, without referring to the underlying scaling law of geographic features.

Despite the applicability of fractal thought to cartographic practices, the fractal nature of maps and mapping has not been well received in the cartography literature. This can be seen from the fact that some series maps have some arbitrary scaling ratio, for example, from 1:100K to 1:250K, with a scaling ratio of 1/2.5. This arbitrary scaling ratio creates a problem for map generalization. From 1:500K to 1:1M, we simply reduce each of the four 1:500K maps by a quarter (or a half in terms of the side length), and put the reduced four together as one 1:1M map for map generalization. Similarly, from 1:150K to 1:450K, we simply reduce each of the nine 1:150K maps by one ninth (or one third in terms of the side length), and put the resulting reduced nine together as one 1:450K map for generalization. Along this line of logic, we cannot figure out how to generalize a map from 1:100K to 1:250K. Therefore, if we were explicitly guided by fractal geometry, we would not have had the arbitrary scaling ratio in the series maps. This is something we can learn from the fractal nature of maps. Given the fractal nature of maps, mapping can be considered as a head/tail breaks process.

## 4. Mapping as a head/tail breaks process

The so-called head/tail breaks process involves dividing things around an average into large and small, which respectively constitute the head and the tail of the rank-size plot (Zipf 1949), and recursively continuing the dividing process for the large until the notion of far more small things than large ones is violated (Jiang 2013a). Note that the terms 'large' and 'small' should be understood broadly, representing 'popular' and 'unpopular' in terms of topological property, or 'meaningful' and 'meaningless' in terms of semantic property. Depending on the scales (or sizes), things can be classified or clustered into different hierarchical levels. However, unlike conventional classifications, head/tail breaks is able to reveal the underlying fractal or scaling structure, thereby being unique and powerful for mapping.

### 4.1 Map generalization guided by head/tail breaks

We will now examine how map generalization can be conducted following head/tail breaks. A 'great wall' looks like the third iteration Koch curve on the 1:50K map (Figure 3). The great wall or Koch curve contains three discrete scales in a recursive manner: 1/3, 1/9, and 1/27, where the scaling ratio is 1/3. In order to generalize the wall or curve and represent it in the 1:150K and 1:450K maps, we first reduce the 1:50K map by one third, and then select the large line segments (1/3, and 1/9) in order to get a map of 1:150K. To further represent the curve in the 1:450K map, we reduce the 1:50K by one ninth (the left panel of Figure 3), or the 1:150K curve by one third (the right panel of Figure 3), and then select the largest line segments of 1/3. In essence, the process is based on head/tail breaks; in other words, by recursively selecting large things in the head, or equivalently eliminating small things in the tail. The reader might have noted that the generalized curves are in fact the reduced curves of the second and first iteration, as is clearly illustrated in Figure 3.

We will now further examine two patterns, known as the Sierpinski carpet and the Sierpinski triangle, from the source map of 1:50K. The scaling ratios for the two patterns are 1/3, and 1/2, respectively. Figure 4 illustrates the generalized results for the two series maps: 1:50K, 1:150K, and 1:450K with the scaling ratio of 1/3, and the other 1:50K, 1:100K, and 1:200K with the scaling ratio of 1/2. The scaling ratios of the two series maps match those of the two patterns; otherwise, the results could not be considered correct, strictly speaking. This is an open issue that will be further discussed in the next section; in other words, how strictly should scaling ratio be respected in the course of map generalization? The generalization looks fairly simple, since both the Koch curve and the Sierpinski patterns involve only three discrete scales. However, despite this simplicity, the examples capture the essence of generalization, i.e., keeping large things while eliminating small things, or equivalently keeping well-connected, popular or meaningful things, while eliminating less-connected, unpopular or meaningless things. More examples of map generalization based on head/tail breaks can be found in Jiang et al. (2013).

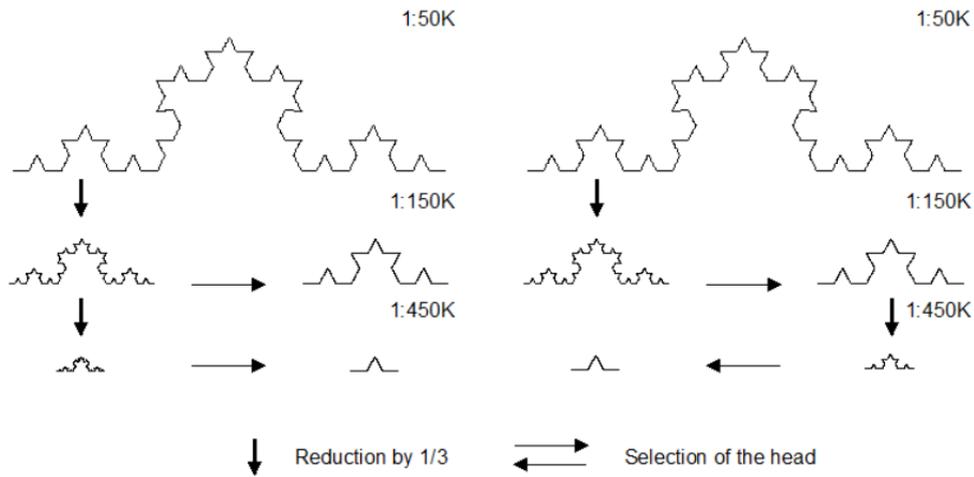

Figure 3: An illustration of map generalization by head/tail breaks
(Note: There is a slight difference between the left and right panels. For the left panel, both 1:150K and 1:450K are derived from the same source 1:50K, whereas for the left panel, 1:150K from 1:50K, and subsequently 1:450K from 1:150K)

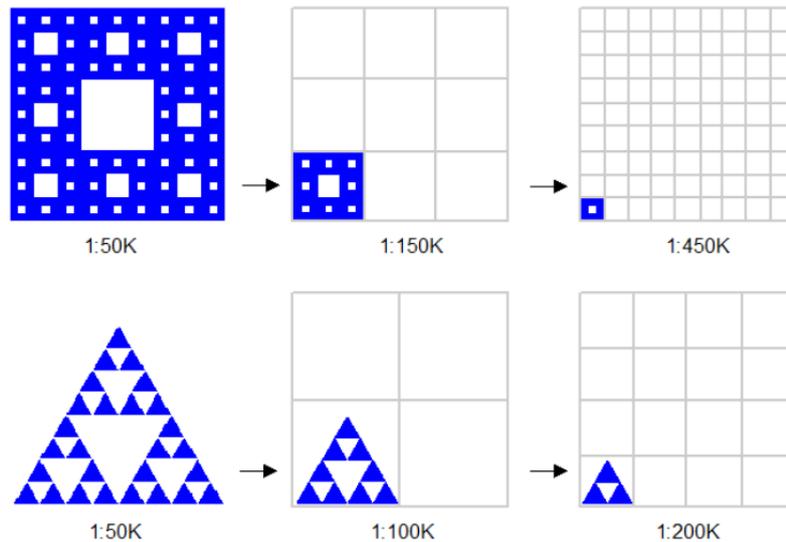

Figure 4: (Color online) Generalization of two patterns: the Sierpinski carpet and Sierpinski triangle
(Note: The background grids indicate scaling ratios of the series maps.)

**4.2 Statistical mapping based on head/tail breaks**
Statistical mapping can benefit considerably from fractal thought, and head/tail breaks in particular. To illustrate this, we assume a set of 1023 cities whose sizes follow Zipf's law exactly, i.e., 1, 1/2, 1/3, 1/4,…, 1/1023 (Figure 5, Panel a), implying that the first largest city is twice as big as the second largest, three times as big as the third largest, and so on. For all of the 1023 cities, their average size is 0.0073, which splits the cities into two parts: those above the mean 137 cities in the head, and those below the mean 886 cities in the tail. The 137 cities in the head can be further split into two parts around the second mean of 0.0402: this yields 24 cities above the mean in the head, and 113 cities below the mean in the tail. This process can recursively continue (see the result in Figure 5, Panel c). Eventually, all of the 1023 cities can be classified into four hierarchical levels, which implies that the pattern of far more small cities than large ones recurs three times (see the nested rank-size plots in Figure 5, Panel b). The scaling pattern revealed is a fairly accurate reflection of the underlying scaling structure. Note that the resulting head percentage is always a minority, and this would not be the case if the conventional classification method natural breaks (Jenks 1967) were used. It is important to observe that the pattern of far more small things than large one recurs again and again, which indicates a sort of self-similarity from the perspective of the nested rank-size plots.

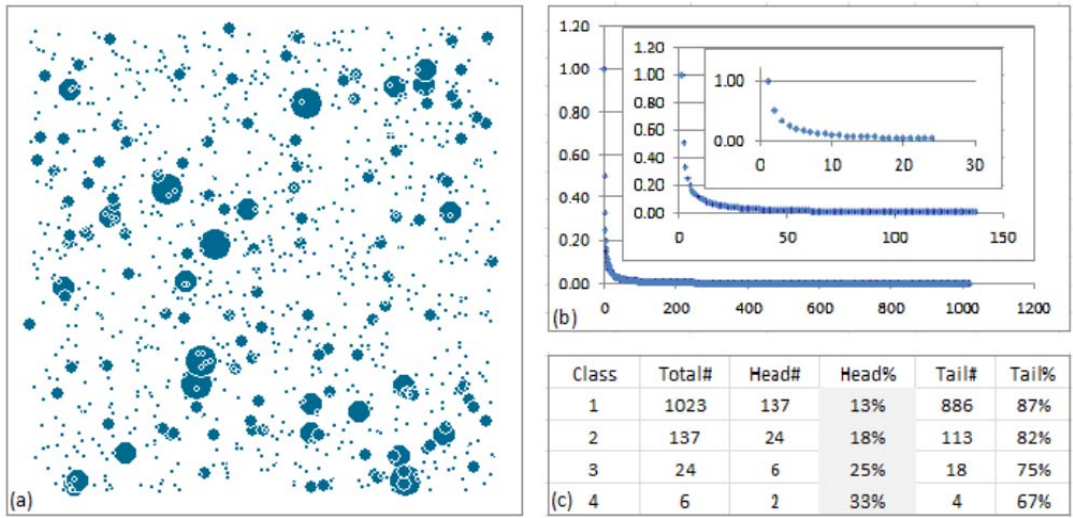

Figure 5: (Color online) The fractal structure of the 1023 cities that follow Zipf's law exactly: (a) the pattern, (b) the nested rank-size plots, and (c) the head/tail breaks statistics
(Note: The city sizes are 1, 1/2, 1/3, 1/4,…, and 1/1023, implying that the first largest city is twice as big as the second largest, three time as big as the third largest, and so on.)

**4.3 Cognitive mapping as an extension of map-making**
Cognitive mapping can be considered an extension of ordinary map-making, and again is governed by head/tail breaks. Given the city pattern shown in Figure 5, what is the corresponding cognitive map in the human mind? There is little doubt that the largest cities, or the six largest cities to be more specific, constitute the major spots of a cognitive map. This cognitive mapping is primarily derived from a geometric perspective, i.e., relies on the size for the head/tail breaks process of geographic features. In many cases, topology and semantics play more important roles than geometry in cognitive mapping (see Appendix). For example, a city may not be the largest, but rather the most-connected among a city system; similarly a city may be neither the largest nor the most popular, but carries the highest semantic meaning, such as being the oldest city. Essentially, the largest, the most popular, or the most meaningful constitutes a cognitive map of cities. Following the classic works on cognitive maps (Tolman 1948, Lynch 1960), a large body of literature has focused on animal or human internal cognitive processes. Instead of seeking reasons internally, Jiang (2013b) argued that it is the underlying fractal structure that makes cities legible or imageable, and further claimed that the fractal structure, or the external representation itself, is the first and foremost cause, without which mental maps would not be easily formed. Fractal structure forms the essence of cognitive mapping, as well as that of cartographic mapping as discussed earlier.

This section, for the sake of simplicity, adopts classic (or strict) fractals to illustrate various mappings, but the underlying spirit applies to statistical fractals, such as coastlines or any other geographic features (Jiang et al. 2013). In this regard, head/tail breaks provides a powerful means by which to derive the inherent scaling hierarchy of geographic features or data for mapping purposes. The use of head/tail breaks may sound rather simplistic, since it does not consider the differences between power laws, lognormal and exponential distributions; head/tail breaks applies as long as there is a small head and a long tail. However, while this simplistic thinking might be of little use from a physicist point of view (e.g., Bak 1996), it is of great value for mapping or understanding the underlying geographic forms and processes. Geographic forms are fractal rather than Euclidean, and the underlying geographic processes are nonlinear rather than linear (e.g., Batty and Longley 1994, Frankhauser 1994, Salingaros 2005, Benguiui and Blumenfeld-Lieberthal 2007, Chen 2009, Jiang and Yin 2014). This presents the basic perception of geographic forms and processes. In the next section, we further discuss the implications of the fractal nature of maps and mapping for cartography and geography.

## 5. Implications of the fractal nature of maps and mapping

A map, or more precisely, a well-designed map, is both a scientific product and an artistic artifact. Following our aforementioned arguments, we can remark that maps that reflect the underlying fractal structure are well-designed, since they are scientifically correct and artistically appealing, with positive impacts on human well-being (Jiang and Sui 2014). However, map-making practices, as well as cartographic research, have not been explicitly guided by fractal geometry, or fractal thought in general. On the contrary, cartographers and geographers alike are misguided by Euclidean geometry, because they tend to (1) focus on individual things or scales, but forget the scaling pattern across all scales, and (2) believe more or less similar things, but forget far more small things than large ones. Herewith we further discuss implications of the fractal nature of maps and mapping.

Geographic features or phenomena can be described as individually with scale, and collectively without scale (or scale-free as a synonym of fractal and scaling). This description represents two distinct ways of looking at geographic features: using Euclidean geometry on the one hand, and fractal geometry on the other. One can measure the sizes of things at an individual level, but collectively there is no typical size for characterizing the sizes of things. Instead, there are all sizes of things, or far more small things than large ones. This scale-free or fractal property of geographic features implies that geographic systems are complex, and that conventional Newtonian physics is not a suitable paradigm for understanding their complexity. Rather, we are required to adopt complexity modeling tools, in particular simulations from the bottom up that focus on interactions among individuals, in order to better understand geographic forms and processes (Jiang 2014).

The fractal nature of maps and mapping makes cartography special, and differentiates maps from other graphical or pictorial representations. A map is a model or an approximation of reality, rather than an image or mirror of reality. Figure 6 illustrates these two different views using the same Koch curve. The power of maps lies in the model view rather than the image view. Strictly speaking, satellite images are not maps, for they do not involve mapping processes such as generalization, classification, and symbolization. It is these mapping processes that make cartography and maps unique and special. Generalization and classification are shared by cognitive mapping, because cognitive maps tend to be schematic rather than detailed, and topological or semantic rather than geometric. Despite the pervasiveness of high-resolution satellite images of the earth surface, both topographic and thematic maps remain largely irreplaceable.

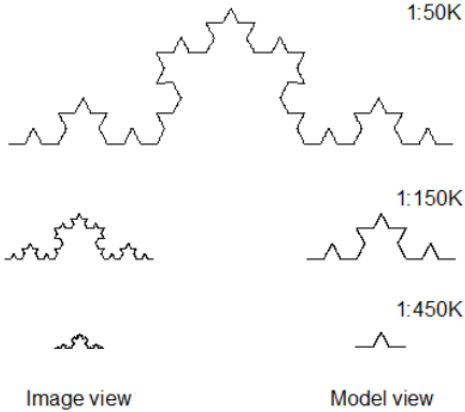

Figure 6: Map as a model (to the right) rather than as an image of reality (to the left)

The fractal nature of maps and mapping implies a holistic view of looking at what is to be mapped, as well as recognition of the underlying fractal structure. The fractal nature also implies moving beyond a geometric perspective towards a topological and/or semantic perspective. For example, a street network can be assessed from a geometric view, in which junctions connect to junctions, or street segments connect to street segments. However, this geometric view, due to the geometric locations of the junctions and distances of the segments, prevents us from seeing the underlying fractal patterns. Instead, a topological view in terms of how individual streets (such as named streets) connect to other streets enables us to see all kinds of streets in terms of connections; in other words, far more less-connected streets than well-connected ones, or equivalently far more unpopular streets than popular ones (Jiang 2007). In parallel with the topological view, a large scope, rather than a small scope, is essential for seeing the underlying fractal structure. For example, a city scope is large enough to see the fractal structure of streets. As mentioned in Section 4.3, semantics plays a more important role than geometry in mapping (see the Appendix for an example). Eventually, mapping involves all geometric, topological and semantic information, and recognizes the underlying scaling property in particular, before applying the head/tail breaks process.

The scaling ratio has some profound implications for mapping, or map generalization, statistical mapping in particular. In series maps subdivision, the ratio of the two subsequent map scales should be a reciprocal of an integer like 1/2 and 1/3 rather than that of some arbitrary numbers like 1/2.5. For example, the scaling ratio of the series maps 1:100K, 1:250K, and 1:625K is 1/2.5. This implies that from 1:100K to 1:250K, or from 1:250K to 1:500K, the number of map sheets is neither from 4 to 1, nor from 9 to 1. This is problematic, and is probably a disadvantage when the series map subdivision is not explicitly guided by fractal geometry. In order to further clarify this point, we re-assess the map generalization of the Sierpinski carpet and Sierpinski triangle shown in Figure 4. Note in this example that the scaling ratios of the Sierpinski carpet and triangle (1/2 and 1/3, respectively) are the same as those of the map scales. How can we generalize the Sierpinski carpet (1:50K) into map scale 1:100K, and 1:200K, or the Sierpinski triangle (1:50K) into map scale 1:150K, and 1:450K? Strictly speaking, the results would not be the same as those shown in Figure 4. We have no answer to this question, and it deserves further research in the future.

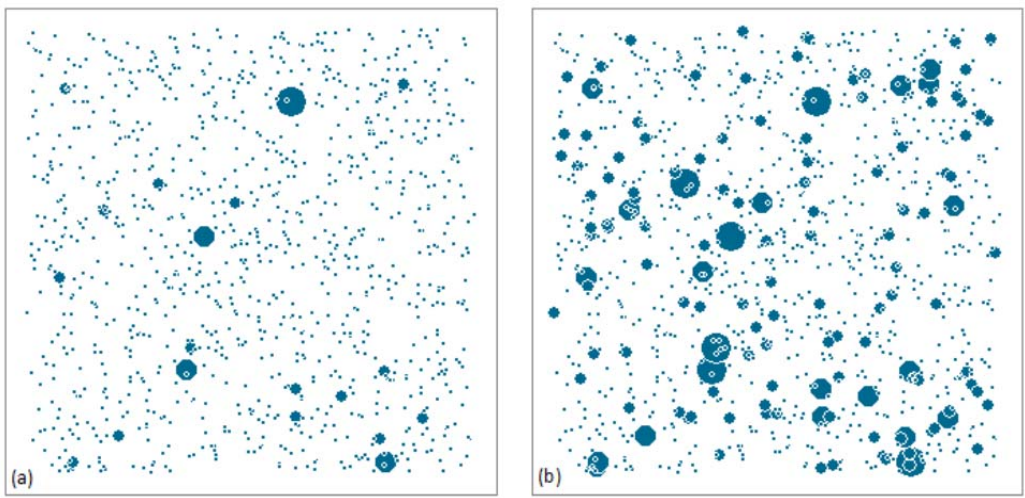

Figure 7: (Color online) The 1023 cities that follow Zipf's law exactly, yet using different classifications of: (a) natural breaks, and (b) head/tail breaks
(Note: The city sizes are 1, 1/2, 1/3, 1/4,…, and 1/1023, implying that the first largest city is twice as big as the second largest, three time as big as the third largest, and so on.)

A map that reveals the underlying fractal structure, i.e., based on head/tail breaks, can evoke a sense of beauty (Jiang and Sui 2014). This new kind of beauty, initially discovered and defined by Christopher Alexander (1993, 2002), lies in the underlying fractal structure, and is therefore related to objective judgments rather than subjective preferences. In order to verify the existence of this beauty, Alexander (1993, 2002) developed what he called the 'mirror of the self' test, using Turkish carpets and building pictures as testing targets. In the test, human subjects were asked to look at two patterns placed side by side and identify which of the two most resembled the subject's inner self, or character. Interestingly, most of the time, and for most of people, it was found that patterns with fine structures (or fractal patterns, as we refer to them in this paper) have a better sense of the 'mirror of the self' feeling. For example, given the two city system patterns based on the data discussed in Section 4.2, which respectively used natural breaks (Jenks 1967) and head/tail breaks (Figure 7), the reader is invited to run a 'mirror of the self' test by asking which pattern resembles your own character. We suggest take a few minutes to think about this before turning to Note 1 at the end of this paper for the answer.

The fractal nature of maps and mapping has some implications for map design, as well as map reading. For example, map symbol sizes, color rendering, and label sizes should all be guided by head/tail breaks, or fractal thinking in general. Alexander's 'mirror of the self' test can be of great value for assessing and comparing map designs. However, due to the space limitation of this paper, we cannot explore the topic in further detail; therefore, it should be further researched in the future. The last point to make before presenting our conclusion is that unlike paper maps created prior to the digital era, maps created during the digital era, or the big era in particular, should be increasingly used for illustrating the underlying scaling patterns rather than for measuring things. This is because measuring can be achieved at the database level for some specific map scale, or simply because geographic features in general are not measurable due to the conundrum of length.

## 6. Conclusion

A map is not the territory it represents, but a simplified version of the territory based on mapping processes such as generalization, classification, and symbolization. This paper argues and provides evidence for the fractal nature of maps and mapping, which is fundamentally the same as that of geographic features, that is, with far more small things than large ones. To remind, this notion of far more small things than large ones is not just constrained to geometric sizes or details, but should be broadly interpreted as far more unpopular things than popular ones, and far more meaningless things than meaningful ones, respectively in terms of topological and semantic properties. Fractal structure is ubiquitous not only in rocks, rivers and watersheds, mountains, islands, and coastlines, but also in human-made artifacts such as cities, streets, buildings, social media, the Internet, and the World Wide Web. All of these constitute targets to be mapped in current cartography. A map is a model, rather than an image of reality, and mapping processes must be guided by fractal geometry, or the notion of far more small things than large ones. Fractal geometry also provides some theoretical backing for cartographic principles such as series maps subdivision, visual hierarchy, and Töpfer's radical law.

Cartography is primarily founded on Euclidean geometry in terms of representing geographic features using geometric primitives such as points, lines, and polygons. It is also Euclidean geometry, focusing on individual scales, that prevents us from seeing the underlying scaling pattern across all scales. Shifting from geometry to topology or from geometry to semantics helps us to see the underlying fractal structure, or the recurring scaling pattern of far more small things than large ones. Maps, to a great extent, are to reveal the underlying scaling structure. In this regard, mapping can be considered the head/tail breaks process for dividing things into different hierarchical levels (for statistical mapping), into large and small things (for map generalization), and for identifying the largest things (for cognitive mapping). Maps and fractals share common properties such as recursion, self-similarity, constant scaling ratio, and characteristic scaling exponent. Scaling must be formulated as a law of cartography or that of geography in general.

Note 1: According to Alexander's theory of centers, the right one would have a better sense of the 'mirror of the self' feeling.

**Acknowledgement**
I would like to thank the anonymous referees and the editor Brian Lees for useful comments. I am also indebted to Nikos Salingaros, Manfred Buchroithner, Gilles Benguigui, Georg Gartner, Tony Moore, and Ling Bian who, while not necessarily agreeing on all that I have argued here, provided helpful suggestions for revisions.

**Appendix: Mapping the fractal structure of this paper**

There are far more low-frequent words than high-frequent words in any text. More precisely, given any naturally generated text such as a book or a paper, the frequency of any word is inversely proportional to its rank, with the most frequent word at number one, the second most frequent word at number two, and so on. Put simply, the most frequent word occurs twice as often as the second most frequent word, three times as often as the third most frequent word, and so on. This statistical regularity, together with city size distribution, was formulated as Zipf's law. This appendix introduces a word cloud by which to map the underlying fractal structure of this paper.

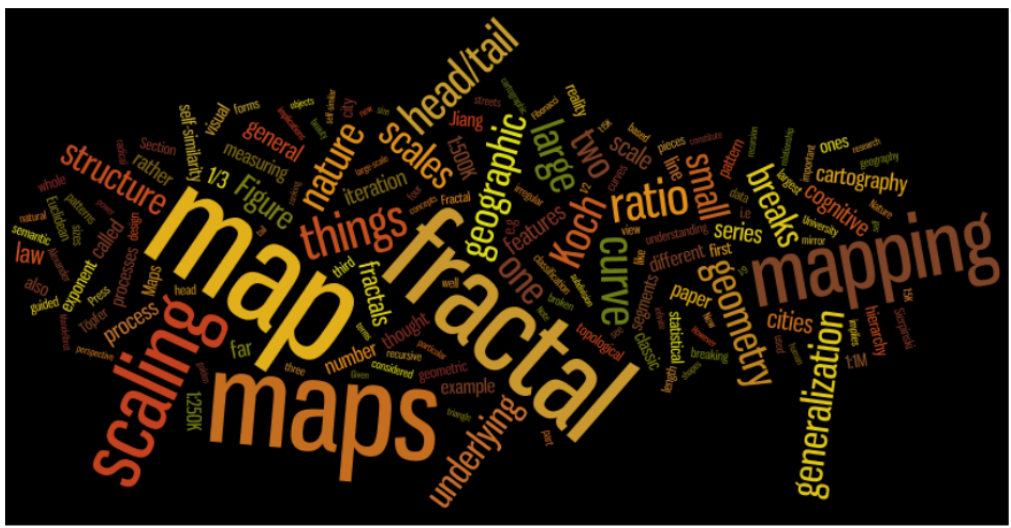

Figure A1: (Color online) The map that captures the underlying fractal structure of the paper

This paper (excluding the appendix) contains 7389 words, but only a few words recur very frequently. The word cloud or map (Figure A1) was created based on the text in this paper using Wordle

(http://www.wordle.net/). In this map, the most frequently used words are shown as map (or maps), fractal, scaling, and mapping; there is no doubt that the word frequency follows Zipf's law, or a scaling hierarchy. The frequently recurring words are in the head, while the many others are in the tail of a heavy tailed distribution in a rank-size plot. Note that this map is neither geometric nor topological, but semantic. This is because in the map the locations of the words, and the distances or directions between the words, are all arbitrarily determined. No topological relationship between the words is mapped. Instead, the word sizes reflect some semantic meaning in terms of how they are related to the thesis of this paper; the most meaningful words have the largest sizes.

The fractal structure of the paper can be further seen at the different levels of scale; that is, the paper is composed of six sections, which are further decomposed into 31 paragraphs, 252 sentences, and 7389 words, which display a striking scaling hierarchy. There are three intermediate scales (sentences, paragraphs, and sections) between the smallest (words) and the largest (the paper) scales. In fact, the words can be further decomposed into letters, but we consider words to be the building blocks of the paper. The top four words (maps, mapping, fractal, and scaling) capture this paper's central thesis, which is reflected in the paper's title *'The fractal nature of maps and mapping'*. The map in Figure A1 looks beautiful and eye-catching. However, to paraphrase Christopher Alexander (2002), it is essentially the underlying fractal structure, rather than the colors, that makes the map beautiful.